\documentclass{amsart}
\usepackage{amssymb}
\usepackage[mathscr]{eucal}
\usepackage{bbm}

\usepackage{hyperref}

\allowdisplaybreaks[1]

\newtheorem{thm}{Theorem}

\theoremstyle{definition}

\numberwithin{equation}{section}

\newcommand{\abs}[1]{\left| #1 \right|}
\newcommand{\pr}[1]{\left( #1 \right)}

\newcommand{\R}{\mathbb{R}}

\newcommand{\N}{\mathbb{N}}

\newcommand{\Cc}{\mathscr{C}}

\thanks{The authors wish to thank Joseph Vandehey for helpful discussions about the hotspot lemma.}

\author[D. Airey]{Dylan Airey}
\address[D. Airey]{Department of Mathematics, Princeton University, Fine Hall, Washington
Road, Princeton, NJ 08544-1000, USA}
\email{dairey@math.princeton.edu}

\author[B. Mance]{Bill Mance}
\address[B. Mance]{Uniwersytet im. Adama Mickiewicza w Poznaniu,
  Collegium Mathematicum, ul. Umultowska 87, 61-614 Pozna\'{n}, Poland}
\email{william.mance@amu.edu.pl}

\begin{document}

\title{Hotspot lemmas for non-compact spaces}

\begin{abstract}
We extend the hotspot lemma, or Piatetski\u{i}-Shapiro normality criterion, to dynamical systems on non-compact spaces. We build on the work of N. G. Moshchevitin and I. D. Shkredov, noting an issue in some of their proofs, and adding a necessary and sufficient condition to remove this issue.
\end{abstract}

\maketitle

The hotspot lemma, first established for the $b$-ary shift on the unit circle in \cite{PS}, is a useful tool for showing the normality of digit expansions. This has been extended to general compact ergodic dynamical systems in \cite{Post} and finer statements for systems such as finite state Markov chains have been established in \cite{MS}. In extending the hotspot lemma to systems with infinitely many digits, such as the continued fraction expansion and infinite Bernoulli shifts, the non-compactness of the symbolic shifts leads to an issue with the proofs of some theorems in \cite{MS}. We add a tightness condition which removes this issue and state corrected versions of the theorems.

Let $X$ be a set with a collection of subsets $\Cc = \{C_m\}$ which form a \textit{semi-$\sigma$-algebra}: that is $\Cc$ contains $X$ and $\emptyset$, $\Cc$ is closed under finite intersection, and for any $A \in \Cc$ there is a countable disjoint collection of sets $\{C_k\} \subseteq \Cc$ such that $X \setminus A = \bigcup_k C_k$\footnote{We are unaware of a standard term for this structure and have decided on this name in analogy with the definition of semialgebra, which requires a finite disjoint union rather than a countable one. Note that the conditions for being a $\sigma$-algebra are stronger than those for being an algebra, while the conditions for being a semi-$\sigma$-algebra are weaker than those for being a semialgebra.}. Endow $X$ with the topology and Borel $\sigma$-algebra generated by $\Cc$. Let $\mu$ be a probability measure and let $T : X \to X$ be a continuous map which preserves $\mu$ and is ergodic with respect to $\mu$.

The \textit{Birkhoff mean} of a measurable function $f$ with respect to a point $x \in X$ is given by
$$
S_N(x_0,f) = \frac{1}{N} \sum_{n=0}^{N-1} f(T^i x_0).
$$
We define the sets
$$
A_\ell(f,\delta) = \{x \in X : \abs{S_\ell(x,f) - \int f d\mu} > \delta \}.
$$

We define an analog of the Hausdorff measure $H( \cdot )$ for a set $E$ with respect to this family to be $\inf\{\sum \mu(C_i)\}$,
where the infimum is taken over coverings (finite or countable) of $E$. We say that the measures $\mu$ and $H$
are {\it coordinated} if any $\mu$-measurable set is $H$-measurable.

Theorems $1$, $4$, and $5$ in \cite{MS} are incorrect as stated. A counterexample to all three is given by considering the space $X = \N^\N$ with the shift map $T$ and the family $\Cc$ given by the cylinder sets $[\xi] = \{x \in X : x|_{\{1, \cdots, |\xi|\}} = \xi\}$ for $\xi \in \N^{<\infty}$. Consider the point $x_0 = (1,2,3,4,\cdots)$. Then for any $I \in \Cc$
\begin{align*}
\limsup_{N \to \infty} \frac{S_N(x_0, \chi_I)}{N} = 0
\end{align*}
since for a fixed $\xi \in \N^{<\infty}$, if $M = \max_{1 \leq i \leq |\xi|} \xi_i$, then for $n > M$ we have $T^n \notin [\xi]$. For any probability measure $\mu$ on $X$ and function $\varphi : \R_{\geq 0} \to \R_{\geq 0}$ we have $0 \leq \varphi(\mu(I))$. Thus $x_0$ satisfies the assumptions in each theorem. However, for any probability measure $\mu$ on $X$, there must be some $a \in \N$ such that $\mu[a] > 0$. Thus,
\begin{align*}
\lim_{N \to \infty} \frac{S_N(x_0,\chi_{[a]})}{N} \neq \mu[a].
\end{align*}

The error in the proofs appear in equations $(6)$ and $(19)$ for Theorems $1$ and $4$ respectively, and Theorem $5$ uses Theorem $1$. Specifically, the authors distribute a $\limsup$ over an infinite sum but in general
$$
\limsup_n \sum_i a_{i,n} \leq \sum_i \limsup_n a_{i,n}
$$
only for finite sums. For Theorems $2$ and $3$ the spaces in consideration are compact, so these sums will be finite while for Theorem $1$ the space is not necessarily compact and for Theorems $4$ and $5$ the spaces are non-compact.

The behavior of this counterexample where the mass of the orbit escapes to infinity is the only obstruction to the theorems. To correct these theorems we add a tightness condition which prevents this escape. We say a set of probability measures $M$ is \textit{tight} if for every $\epsilon > 0$ there is a compact set $K$ such that for every $\mu \in M$, $\mu(X\setminus K) < \epsilon$. We define the empirical probability measures for $x \in X$ by $\mathcal{E}(x,n) = \frac{1}{n} \sum_{i=0}^{n-1} \delta_{T^i x}$. Note that $S_N(x_0, f) = \int f d\mathcal{E}(x_0,N)$.
\begin{thm}[Correction of Theorem~$1$ in \cite{MS}]
Let $x_0 \in X$ be such that the set of probability measures $\{\mathcal{E}(x,n)\}_{n=1}^\infty$ is tight and let $\varphi : R_{\geq 0} \to R_{\geq 0}$ be monotone increasing. If for an arbitrary set $I$ from the family $\Cc$
\begin{align*}
\limsup_{N \to \infty} \frac{S_N(x_0, \chi_I)}{N} \leq \varphi(\mu(I))
\end{align*}
and for any $\delta > 0$ we have $\lim_{\ell \to \infty} H_\varphi(A_\ell(\chi_I,\delta)) = 0$, then for a Borel set $B$, we have
\begin{align*}
\lim_{N \to \infty} \frac{S_N(x_0, \chi_B)}{N} = \mu(B).
\end{align*}
\end{thm}
Note that when $X$ is compact, the set of empirical measures $\{\mathcal{E}(x,n)\}_{n=1}^\infty$ is tight for every $x \in X$, so Theorem $1$ in \cite{MS} is true for compact spaces.
While we will record the proof in full for this theorem, the modifications to the proof of Theorem $4$ are identical. Thus, we omit them.
\begin{proof}
Let $\delta > 0$ and let $I$ be a set in $\Cc$. Pick a compact set $K \subseteq X$ such that $\mathcal{E}(x_0,n)(X \setminus K) < \delta$ for every $n \in \N$. For all positive integers $N$ and $\ell$ we have
\begin{align*}
S_N(x_0, \chi_I) = \frac{1}{\ell} \sum_{i=0}^{N-1} S_\ell(T^i x_0, \chi_I) + O(\ell).
\end{align*}
Thus,
\begin{align*}
\abs{\frac{S_N(x_0, \chi_I)}{N} - \mu(I)} &= \frac{1}{N} \abs{\sum_{i=0}^{N-1} \pr{\frac{S_\ell(T^i x_0, \chi_I)}{\ell} - \mu(I)}} + O\pr{\frac{\ell}{N}}
\end{align*}
We decompose this sum into three parts: first the set of indices $i$ such that
\begin{align*}
\abs{\frac{S_\ell(T^i x_0, \chi_I)}{\ell} - \mu(I)} < \delta,
\end{align*}
second the set of indices $i$ such that
\begin{align*}
T^{i+j} x_0 \notin K \text{ for some } 0 \leq j \leq \ell -1,
\end{align*}
and third the set of indices $i$ such that
\begin{align*}
T^{i+j} x_0 \in K \text{ for all } 0 \leq j \leq \ell-1 \text{ and } \abs{\frac{S_\ell(T^i x_0, \chi_I)}{\ell} - \mu(I)} \geq \delta.
\end{align*}
We have $
\abs{\frac{S_\ell(T^i x_0, \chi_I)}{\ell} - \mu(I)} \leq 2$
which implies
\begin{align*}
\abs{\frac{S_N(x_0, \chi_I)}{N} - \mu(I)} \leq \delta + \frac{2 Q_N}{N} + \frac{2 R_N}{N} + O\pr{\frac{\ell}{N}}
\end{align*}
where $Q_N$ is the number of terms in the second sum and $R_N$ is the number of terms in the third sum. Now
\begin{align*}
Q_N &\leq \sum_{j = 0}^{\ell -1} \#\{0 \leq i \leq N-1 : T^{i+j} x_0 \in X \setminus K\} \\
& \leq \ell \#\{0 \leq i \leq N+\ell-1 : T^i x_0 \in X \setminus K\} \\
& = \ell (N+\ell -1) \mathcal{E}(x_0, N+\ell -1)(X\setminus K) \\
&< \ell (N+\ell-1) \delta.
\end{align*}
For every index $i$ in the third sum we have $T^i x_0 \in A_\ell(\chi_I,\delta) \cap K$. Choose a countable collection of sets $\{C_k\}$ from $\Cc$ which cover $A_\ell(\chi_I, \delta) \cap K$ and such that
$$
\sum_i \varphi(\mu(C_i)) < H_\varphi(A_\ell(\chi_I,\delta) \cap K) + \delta.
$$
Note that every set in $\Cc$ is clopen: each $I \in \Cc$ is open by the choice of topology, and the complement of $I$ is the disjoint union of elements in $\Cc$ and so is open itself. Defining $D_0 = I$ and $D_1 = X \setminus I$, we can decompose the function
$$
S_N(x,\chi_I) = \sum_{\sigma \in \{0,1\}^N} |\sigma^{-1}(0)|\chi_{D_{\sigma(0)} \cap T^{-1} D_{\sigma(1)} \cap \cdots \cap T^{-(N-1)} D_{\sigma(N-1)}}(x).
$$
Both $D_0$ and $D_1$ are clopen and $T^{-i} D_j$ is clopen for any $i \in \N$ and $j = 0,1$ since $T$ is continuous. Thus $S_N(x,\chi_I)$ is the weighted sum of indicator functions of disjoint clopen sets and therefore $A_\ell(\chi_I,\delta)$ is a finite union of clopen sets and is clopen as well. Therefore $A_\ell(\chi_I,\delta) \cap K$ is compact and we may extract a finite subcover $\{C_{k_i}\}$ from $\{C_k\}$. 
By assumption, $H_{\varphi}(A_\ell(\chi_I,\delta) \cap K) \leq H_\varphi(A_\ell(\chi_I,\delta)) \to 0$ as $\ell \to \infty$. Thus,
\begin{align*}
\frac{R_N}{N} \leq  \sum_i \frac{S_N(x,\chi_{C_{k_i}})}{N} \leq \sum_i \varphi(\mu(C_{k_i})) \leq \sum_k \varphi(\mu(C_k)) < 2\delta
\end{align*}
for sufficiently large $\ell$. Taking $N$ sufficiently large so that $\ell/N <\delta$ we have $\abs{\frac{S_N(x_0, \chi_I)}{N} - \mu(I)} < \delta + 2\delta + 4\delta + \delta = 8\delta$. Since $\delta > 0$ was arbitrary we have established the claim for sets $I \in \Cc$. The equality for arbitrary Borel sets then follows from standard approximation arguments.

\end{proof}

\begin{thm}[Correction of Theorem 4 in \cite{MS}]
Let $T$ be the continued fraction map on $X = [0,1)$ with the Gauss measure $\mu(A) = \frac{1}{\log 2} \int_A \frac{1}{1+x}dx$. Let $x \in [0,1)$ be such that $\{\mathcal{E}(x,n)\}_{n=0}^\infty$ is tight. Let $\psi : \R_{>0} \to \R_{\geq 0}$ satisfy $\psi(t) = O(e^{\eta \sqrt{\log 1/t}})$ as $t\to 0$ for any $\eta > 0$. If for any cylinder set $I$
$$
\limsup_{N \to \infty} \frac{S_N(x,\chi_I)}{N} \leq \mu(I) \psi(\mu(I))
$$
then for any Borel set $B$
$$
\lim_{N \to \infty} \frac{S_N(x,\chi_I)}{N} = \mu(I).
$$
\end{thm}

\begin{thm}[Correction of Theorem 5 in \cite{MS}]
Let $p = (p_a)_{a=1}^\infty$ be a probability vector. Suppose for some $\eta_0 > 0$ the series
$$
\sum_{a=1}^\infty p_a^{1-\eta_0}
$$
converges. Consider the system $(X,T,\mu)$ where $X = \N^\N$, $T$ is the right shift, and $\mu$ is the Bernoulli measure given by sampling each digit i.i.d. according to $p$. That is $\mu[a_1, \cdots, a_n] = \prod_{i=1}^n p_{a_i}$. Let $\varphi : \R_{\geq 0} \to \R_{\geq 0}$ be a function such that $\varphi(t) = O(t^{1-\eta})$ as $t \to 0$ for some $\eta \in (0,1)$. Suppose $x \in X$ is such that $\{\mathcal{E}(x,n)\}_{n=0}^\infty$ is tight. If for any cylinder set $I$
$$
\limsup_{N \to \infty} \frac{S_N(x,\chi_I)}{N} \leq \varphi(\mu(I)),
$$
then for any Borel set $B$
$$
\lim_{N \to \infty} \frac{S_N(x,\chi_B)}{N} = \mu(B).
$$
\end{thm}

\end{document}